\documentclass[12pt]{amsart}

\usepackage{amsmath,amssymb}

\begin{document}

\newtheorem{lem}{Lemma}
\newtheorem{cor}{Corollary}
\newtheorem{pr}{Proposition}
\newtheorem{thm}{Theorem}

\theoremstyle{definition}
\newtheorem{ex}{Example}
\newtheorem{rem}{Remark}

\author{S. P. Novikov}
\title[Dynamical Systems and Differential Forms]{\hbox{Dynamical Systems and Differential Forms.}
\vskip3pt
Low Dimensional Hamiltonian Systems}
\address{University of Maryland, College Park, USA;\hfill\break
\phantom{La}Landau Institute, Moscow, Russia,\hfill\break
\phantom{La}fax  301 3149363, tel 301 4054836(o),\hfill\break
\phantom{La}{\it Homepage\/}: {\tt www.mi.ras.ru/\~{}snovikov}}
\email{novikov@ipst.umd.edu}

\maketitle

\centerline{\em Dedicated to Misha Brin}

\vspace{1cm}

\noindent
 {\bf Abstract:} {\it The theory of differential forms began with a discovery of Poincar\'e
who found conservation laws of a new type for Hamiltonian
systems\,---\,The Integral Invariants. Even in the absence of
non-trivial integrals of motion, there exist invariant differential
forms: a symplectic two-form, or a contact one-form for geodesic
flows. Some invariant forms can be naturally considered as ``forms
on the quotient.'' As a space, this quotient may be very bad in the
conventional topological sense. These considerations lead to an
analog of the de Rham cohomology theory for manifolds carrying
smooth dynamical system.  The cohomology theory for quotients,
called ``basic cohomology'' in the literature, appears naturally in
our approach. We define also new exotic cohomology groups associated
with the so-called cohomological equation in dynamical systems and
find exact sequences connecting them with the cohomology of
quotients. Explicit computations are performed for geodesic and
horocycle flows of compact surfaces of constant negative curvature.
Are these famous systems Hamiltonian for a 3D manifold with a
Poisson structure? Below, we discuss exotic Poisson structures on
3-manifolds having complicated Anosov-type Casimir foliations. We
prove that horocycle flows are Hamiltonian for such exotic
structures. The geodesic flow is non-Hamiltonian in the 3D sense.}

\vspace{0.3cm}

\section{Dynamical Systems and Differential Forms. Exact sequences}

\vspace{0.3cm}

In our previous work \cite{N}, we studied various metric independent cohomology groups defined
by subcomplexes in the de Rham complex of differential forms $\Lambda^*(M)$ with differentials
of the form $d_A=d+A$ for a $0$-order operator $A$ acting on differential forms.
In particular, the following two examples were treated:

1. The operator  $d_{\omega}$ defined by
$$d_{\omega}(u)=du+\omega\wedge u,$$ or a family of such
differentials $d_{\lambda\omega},\lambda\in{\mathbb C}$. Such
families were first considered in 1986 in order to find a correct
${\mathbb Z}_2$-graded analog of Morse inequalities for vector
fields extending the ``fermionic technic'' invented by Witten (see
the references in \cite{N}). In this ``metric dependent version,''
we considered the behavior of the zero modes
$$b^{\pm}_{\infty}=\max_{g_{ij}}\limsup_{|\lambda|\rightarrow\infty}b^{\pm}(\lambda),\; \lambda\in\mathbb R,$$
for the family of operators $(d_{\lambda\omega}+d_{\lambda\omega}^*)^2$ acting on the spaces of even
and odd differential forms $\Lambda^{\pm}(M)$.
Here the form $\omega$ corresponds to a vector field $X$ via $\omega_i=g_{ij}X^j$.
We proved that $$b^{\pm}_{\infty}\leq X^{\pm}$$ where $X^{\pm}$ are the numbers
of critical points $X=0$ with signs~$\pm$.

In the metric independent version, the subcomplex $\Lambda_{\Omega}\subset \Lambda^*(M)$
is defined  by the equations $$\Omega\wedge u=0,\quad \Omega=d\omega$$
for all $\lambda\in\mathbb C$. For $\lambda=0$, we have the usual operator $d$
restricted to the smaller subcomplex $\Lambda_{\Omega}\subset \Lambda^*(M)$.
This complex appears naturally on any non-degenerate energy level $H=const$
in a (non-compact) symplectic manifold $(M,\Omega)$ with exact symplectic form $\Omega=d\omega$
(and on all contact manifolds $(M',\omega)$ as well). In particular, on a 3-manifold
$M_c\subset M$ with $\Omega_c=\Omega|_{M_c}$ one gets
$$\Lambda_{\Omega}^0=0,\qquad \Lambda^j_{\Omega}=\Lambda^j(M), j=2,3,$$
and  $u\in\Lambda^1_{\Omega}\subset \Lambda^1(M)$ if and only if $(u,X)=0$.
Here $X$ is the Reeb vector field for the form $\Omega$ on $M$, and $(u,X)=i_X(u)$
is the scalar product (pairing) of a 1-form and a vector field.

\smallskip
2. Another special case of the operator $A$ considered in \cite{N} is the following.
For a vector field $X$, the operator $i_X:\Lambda^k\rightarrow\Lambda^{k-1}$ is defined
on a differential $k$-form locally written as $v=\sum v_{i_1,i_2...i,_,k}dx^{i_1}\wedge... \wedge dx^{i_k}$
by the ``index-summation'' formula $$(X,v)=i_X(v)_{i_1,..,.i_{k-1}}= X^iv_{i,i_1,...,i_{k-1}}.$$
This operator satisfies $i_X^2=0$ so that $$Im(i_X)\subset Ker(i_X)\subset \Lambda^*(M).$$

The square $(d+i_X)^2=\nabla_X$ of the operator $d+i_X$ is the Lie derivative acting
on differential forms. Its kernel $$(d+i_X)^2(u)=\nabla_X(u)=0$$
is the subspace $\Lambda^*_{inv}\subset \Lambda^*(M)$ of
differential forms invariant under the time shifts of the
dynamical system $dx^i/dt=X^i(x)$. In this case, the homology of the operator $d+A$
is called the $X$-invariant homology. The whole pencil $d+\lambda \nabla_X$
is useful here. For a collection of commuting vector fields $X_1,...,X_p$ there is
a multivariable pencil $d+\sum_i\lambda_i\nabla_{X_i}$ where the $\lambda_i$'s are
independent commuting variables. Such constructions have been used to study
actions of compact abelian groups (tori). The so-called {\bf equivariant homology}
can be defined using this pencil.

Let us define two subspaces $Im(i_X),Ker(i_X)\subset \Lambda^k(M)$ consisting of all $C^{\infty}$ differential
$k$-forms $u$ such that $$i_X(u)=(u,X)=0$$ for both subspaces, and $u=(v,X)$ for $u\in Im(i_X)$

For example, it was shown in \cite{N} that the subcomplex $\Lambda_{\Omega_c}\subset \Lambda^*(M_c)$
for the $3$-manifold $M_c$ discussed above has the form
$$0\rightarrow\Lambda^1_X\rightarrow\Lambda^2(M_c)\rightarrow \Lambda^3(M_c)\rightarrow 0$$
with the differential $d_{\lambda\omega}$. For $\lambda=0$, it is
equal to the usual~$d$.

A similar statement is valid in the case of $2k+1$-dimensional
energy levels $M_c\subset M$ in a $2k+2$-dimensional symplectic manifold.
Instead of the form $\omega_c$ we take the form $\omega'=\Omega_c^{k-1}\wedge\omega$, and instead of
$\Omega_c,\Lambda_{\Omega_c}$ we take
$$\Omega'=\Omega_c^k=d(\omega'),\Lambda_{\Omega'}, d_{\lambda \omega'}.$$
Here we have (see \cite{N})
 $$\Lambda^0_{\Omega'}=0,\; \Lambda^1_{\Omega'}=\Lambda^1_X,\; \Lambda^j_{\Omega'}=
 \Lambda^j(M), j\geq 2,$$
where $X$ is the Reeb vector field on an energy level in $(M,\Omega)$.
The case when $\lambda=0$  (i.e.\ when the standard de Rham operator $d$ acts
on a non-standard  subcomplex) was addressed in \cite{N}.

Let $X$ be a smooth vector field on a manifold $M$, and let
$$\nabla_X^k:\Lambda^k(M)  \rightarrow\Lambda^k(M)$$
be the Lie derivative along this field acting on the space
of $C^{\infty}$-differential $k$-forms.
We have $$\nabla^0_X(f)=X^i\partial_if,\qquad f\in C^{\infty}(M),$$
$$\nabla_X^{k}(u)=d(u,X)+ (du,X)=(d+i_X)^2(u).$$
The same operator acting on the subspace of $u\in \Lambda^k_X(M)$
such that $(u,X)=0$ is denoted by $\nabla^{k,*}_X=\nabla^k_X|_{(u,X)=0}$.
For $k=0$, we always have $(f,X)=0$. So $\nabla^{0,*}_X=\nabla^0_X$.
For $k=n$, we have $\nabla^{n,*}_X=0$.

\smallskip
\noindent
{\bf Definitions:}

1. We say that a $k$-form $u\in \Lambda^k(M)$ lies in the
factor-space $\Lambda^k(M/X)$ if $\nabla_X^k(u)=0$ and $(u,X)=0$.
In the literature, such forms are called basic forms for the orbit foliation.
By definition, any such form $u$ is invariant, that is, $u\in\Lambda^*_{inv}(M)$.

2. We say that a $k$-form $u$  lies in the subspace $\Lambda_X^k$
if  $(u,X)=0$ so that $$\Lambda^*(M/X)=\Lambda^*_{inv}\bigcap\Lambda^*_X.$$
Let us further define the spaces
$$Z^k=Ker(d)\subset \Lambda^k(M) \text{ and }  H^k_X=Z^k/d(\Lambda^{k-1}_X).$$
For $ k=1$, our definition implies that $H^1_X=H^1(M)$ because
$\Lambda^0_X=C^{\infty}(M)$.

3. As it  was mentioned above, the invariant de Rham complex
$\Lambda^*_{inv}(M)$ can be identified with the kernel
$$(d+i_X)^2(\Lambda^*_{inv}(M))=0.$$
The subcomplexes $\Lambda_{inv}$ and $\Lambda(M/X)$
with the standard operator $d$ define {\bf the invariant homology groups (ring)
and the homology groups of the factor-space $M/X$}.
Equivariant homology based on pencils such as $d+\lambda\nabla_X$ can also be studied.
The ``Massey operations type'' spectral sequences for this homology can be
easily defined using
power series expansions at $\lambda=0$
in a way similar to the previous work of the present author
(see the references in \cite{N}). We obviously have $\Lambda^n(M/X)=0$.
(Here $n$ is the dimension of $M$.)

Actually, the factor-space $M/X$ may be very bad. In particular,
its  ring of  $C^{\infty}$-functions $\Lambda^0(M/X)=Ker (\nabla^0_X)\subset C^{\infty}(M)$
contains only constants for an ergodic flow $X$. This is true even for the Hilbert space $L_2(M)$.
At the same time, the de Rham complex $(\Lambda^*(M/X),d)$ may still be nontrivial.

\begin{ex}
For an ergodic straight-line flow $X$ on the torus $M=T^n$ (i.e., having generic fully
irrational angles), the complex $\Lambda^*(M/X)$ consists of constant $k$-forms $u$
such that $(u,X)=0$. The invariant complex consists of all constant forms $u\in\Lambda^*(M)$.
The dimensions of the homology groups are
$$\dim(H^k(T^n/X))=(n-1)!/k!(n-k-1)!=
\dim (\Lambda^k(T^n/X)))$$ for the factor-space. The entire
cohomology of the torus is invariant, i.\,e.,
$H^*(T^n)=H^*_{inv}(T^n)$.
\end{ex}

\begin{ex}
Consider a $2n$-dimensional completely integrable Hamiltonian system
$(M,X)$ with compact common level sets of commuting integrals.
Generically, these sets are $n$-tori with straight-line flows. We
obtain a non-trivial ring of differential forms on the (usual
topological) $n$-dimensional factor-space whose points correspond to
the orbits of all the commuting flows. However, the ring
$\Lambda^*(M/X)$ is much bigger. Its elements look like forms on the
usual topological factor-space with values in the sheaf whose fibres
are the spaces $H^*(T^n/X)$ for the generic levels of commuting
integrals. These levels are generically Liouville tori $T^n$
considered in the previous example. However, degenerations of these
tori may seriously affect the structure of that sheaf. There exists
a special case in which the ``action variables'' (canonically
adjoint to the ``angles'') are  globally well-defined smooth
functions $S_1,\ldots,S_n$ on the symplectic manifold $M$. This case
was frequently considered by geometers. Here we get a global
Hamiltonian action of the compact abelian group $G=T^n$. This
special case looks easier. I believe that one can extract all
calculations for this case from the studies of symplectic geometers
of the last two decades. However, this literature has never
considered anything which cannot be reduced to the actions of
compact groups, even completely integrable systems with compact
generic tori.
\end{ex}

\begin{ex}
Let $M_c=(H=c)\subset M$ be an energy level as above, and $X$ the Reeb vector field
(i.e., the Hamiltonian flow on $M_c$). We have $\Omega\in \Lambda_{inv}(M)$,
$i_X(\Omega_c)=(\Omega_c,X)=0$ on $M_c$,
and $\nabla_X(\Omega)=0$. So the forms $\Omega_c^j$ lie in $\Lambda^*(M_c/X)$, and $\Omega$
lies in $\Lambda_{inv}(M)$.

For the special case of a geodesic flow, we have
$$M=T_*(Q),\; H=\sqrt{g^{ij}p_ip_j},\; M_c=(H=c>0),\; \Omega=\sum dp_i\wedge
dx^i.$$
The 1-form $\omega=\sum p_idx^i$ is $X$-invariant, that is, $\omega\in \Lambda_{inv}(M)$.
All forms $\Omega_c^j,j=0,1,\ldots,n-1,$ are in $\Lambda^*(M_c/X)$.
This result is valid for all Hamiltonian functions  $H(p,x)$ of Maupertuis--Fermat--Jacobi type
(i.e., such that $H(sp,x)=sH(p,x), s>0$, as in the case of a Riemannian or Finsler length,
maybe with an additional magnetic term).
 \end{ex}

To justify our definitions above, we need the following easy lemma.

\begin{lem}
{\rm 1.\ }For any $C^{\infty}$-manifold $M$ and a smooth vector field $X$ on it,
the subspaces of differential forms $\Lambda^*_{inv}(M),\;\Lambda^*(M/X)\subset \Lambda^*(M)$
are closed under the action of the de Rham operator~$d$.
{\rm 2.\ }The subspace $\Lambda_X$  is closed under the action of the Lie derivative
$\nabla_X^k:\Lambda^k\rightarrow\Lambda^k$ (we denote this operator acting on $\Lambda_X$ by $\nabla^*_X)$.
The image of $\nabla^{k,*}_X$ is contained in $Im(i_X)\subset \Lambda^{k,*}_X$.
\end{lem}

\noindent
{\it Proof.}
For any $u\in \Lambda^*_X$, we have $\nabla^k_X(u)=d(u,X)+ (du,X)$ where $(u,X)=0$
and $(du,X)=i_X(du)\in Im(i_X)\subset Ker (i_X)$. This proves the second statement.
For the factor-space we have $(u,X)=0$ and $\nabla_X(u)=d(u,X)+ (du,X)=0$ by definition.
So we conclude that $(du,X)=0$. Our statement is obvious for invariant forms.
We also have $\nabla^k(u)=(du,X)+d(u,X)=i_X(du)$ if $(u,X)=0$. The lemma is proved.

\smallskip
As was explained to the present author by Sabir Gusein-Zade (see \cite{EGS}),
we can essentially use a result from~\cite{G}. For any Riemannian metric,
the operator $i_X$ is dual to the operator of exterior multiplication by the 1-form $*X$
on~$\Lambda^*(M)$. The homology of this operator was studied in~\cite{G}
for a closed 1-form but in fact this restriction is superfluous: all the proofs
work for non-closed 1-forms as well. The factor-groups $Ker(i_X)/Im(i_X)$ vanish
for all $k>0$. For $k=0$, they are equal to $\mathbb R^l$ where $l$ is the sum of
the ``Milnor indices'' of all critical points $X=0$ assuming that these points
are isolated and have finite type (i.e., all Milnor indices are finite).
In the non-degenerate case, the number $l$ is simply equal to the number of critical points.
We assume below that our vector field $X$ is of finite type. It would be nice to extend these
results to the case when the critical points $X=0$ are ``regular enough.''

\smallskip
We are going to study the correct irreducible $C^{\infty}$-analog of
the cohomological equation $\nabla^0_X(f)=g$ for differential forms.
By definition, the {\bf Irreducible Higher Cohomological Equation}
is
$$\nabla^{k,*}_X(u)=v,$$ where $u\in\Lambda_X^k$ and $v\in
i_X(\Lambda^{k+1}(M))\subset \Lambda^k_X$.
Let us also introduce the factor-spaces
$$C^k_X=Im(i_X)/Im(\nabla^{k,*}_X)\subset Coker(\nabla^{k,*}_X)$$
associated with the higher cohomological equation.

By definition, we have $Coker(\nabla^n_X)=\Lambda^n(M)/Im(\nabla_X^n)=
H^n_X$ for $n>1$ because $\nabla^n_X(u)=d(u,X)$, and $Im(i_X)=Ker(i_X)$
in all positive dimensions.

\begin{pr}
The following spaces of $C^{\infty}$ differential forms are isomorphic for any
smooth vector field $X${\rm :}
$$H^n_X=Coker(\nabla^{n}_X)=C^{n-1}_X=Im(i_X)/\nabla^{n-1}_X(\Lambda^{n-1}_X).$$
For a measure-preserving flow on a manifold, we have $H^n_X=Coker(\nabla^0_X)$.
For non-vanishing vector fields (such as geodesic flows on $T_1(N)=M$) we have
$C^k_X=Coker(\nabla^{k,*}_X)$ for all $k\geq 0$ because $Ker(i_X)=Im(i_X)$ in this case.
For generic vector fields $X$ with non-degenerate isolated critical points $X=0$,
the subspaces $C^k_X$ and $Coker(\nabla^{k,*}_X)$ coincide for all $k>0$.
For $k=0$, the subspace $C^0_X$ has finite codimension in $Coker(\nabla^0_X)$
equal to the number of critical points (i.e., $Im(i_X)\subset C^{\infty}(M)$
consists of all functions that vanish at the points where $X=0$).
\end{pr}

The proof of this proposition uses the exact sequences introduced below
and a byproduct of the result of \cite{G} mentioned above. The latter
allows us to calculate the spaces $Ker(i_X)/Im(i_X)$ for vector fields $X$
of finite type. In particular, all these spaces are trivial for non-vanishing
vector fields. For vector fields with generic singularities (and even with
``finite type'' isolated singularities) these spaces are trivial in all
positive dimensions $k>0$. For $k=0$, each critical point gives a contribution
isomorphic to $\mathbb R$ (in the generic case).

In the classical ergodic theory, the spaces $Ker(\nabla^0_X)$ and $Coker (\nabla^0_X)$
have been considered in the Hilbert space $L^2(M)$. Probably the first author who
investigated these operators in specific Sobolev spaces of smooth functions on
Riemann surfaces was G.Forni (see \cite{F}). The space $Ker (\nabla^{0}_X)$ for
a measure-preserving  flow (such as a Hamiltonian one) in the ergodic case
consists of constants only but the space $Coker^0_X$ is very complicated
(see below). In the measure-preserving  case, the isomorphism
$$f\rightarrow f\sigma, f\in C^{\infty}(M), \sigma\in \Lambda_{inv}^n(M)$$
gives an identification $\nabla^0_X=\nabla^n_X$.

Let us emphasize here that we are working with {\it smooth\/} forms and functions
(at least, the image $\nabla^0_X(f)$ should consist of continuous functions and
forms in our constructions). Our theorem is valid for $C^{\infty}$ differential
forms only.

We have already  introduced the important subspaces $C^k_X\subset
Coker(\nabla^k_X)$ entering in our main theorems. They are defined
by setting
$$C^k_X=Im(i_X)/Im(\nabla^{k,*}_X)\subset
Coker(\nabla^{k,*}_X)=\Lambda^k_X/\nabla^k_X(\Lambda^k_X).$$

\begin{thm}
For any smooth vector field $X$ on an $n$-manifold $M$ and every $ k=-1,0,1,\ldots, n-1$,
the following canonical exact sequence is defined:
$$0\rightarrow Z^k(M/X)\rightarrow Ker(\nabla_X^{k,*})\rightarrow Z^{k+1}(M/X)\rightarrow $$
$$\rightarrow H^{k+1}_X \rightarrow C^k_X\rightarrow H^{k+2}_X\rightarrow H^{k+2}(M)\rightarrow 0,$$
where $\nabla_X^k:\Lambda^k_X\rightarrow \Lambda^k_X$,  $(u,X)=0$ for all $u\in \Lambda^k_X$,
and $Ker (\nabla^k_X)=\Lambda^k(M/X)$.
This  exact sequence can also be represented in the following form
$$0\rightarrow H^{k+1}(M/X)\rightarrow  H^{k+1}_X\rightarrow C^k_X\rightarrow
H^{k+2}_X\rightarrow H^{k+2}(M)\rightarrow 0.$$
\end{thm}

\noindent
{\it Proof.}
Let us define the required homomorphisms.

1. Every $k$-form $u\in Z^k(M/X)$ satisfies $du=0$ and $(u,X)=0$ by definition,
so $\nabla^k_X(u)=(du,X)\pm d(u,X)=0$. The definition of an imbedding
$$m_*:Z^k(M/X)\subset \Lambda^k(M/X)= Ker(\nabla^{k,*}_X)$$ is therefore obvious.

2. The operator $d$ defines a map $$d_*:Ker(\nabla^{k,*}_X)\rightarrow Z^{k+1}(M).$$
However,  $du\in Z^{k+1}(M/X)$ because $(du,X)=0$ and $d(u,X)=0$. This follows from
the properties of the Lie derivative already used before.

3. The map $i_*:Z^{k+1}(M/X)\rightarrow H^{k+1}_X$ is just the natural imbedding
of cocycles $$i_*:u\rightarrow u+d(\Lambda^k_X).$$

4. The map $j_*:H^{k+1}_X\rightarrow C^k_X\subset Coker(\nabla_X^{k,*})$ is defined
by setting
 $$j_*:u\rightarrow v=(u,X)=i_X(u)+\nabla^{k,*}_X  (\Lambda^k_X).$$

5. The map $h_*:C^k_X\rightarrow H^{k+2}_X$ is defined by
$$h_*:v=(u,X)\rightarrow du+d(\Lambda^{k+1}_X).$$

6. The last map $g_*:H^{k+2}_X\rightarrow H^{k+2}(M)$ is defined in the trivial natural way.

\smallskip
Let us show now that this sequence is exact.

1. The cocycle property $d_*(u)=0$ obviously implies that $u\in Z^k(M)$.
By using $(u,X)=0$ and $du=0$, we conclude that $u\in Z^k(M/X)$, i.e.,
$du=(u,X)=0$ and $\nabla^k_X(u)=d(u,X)\pm  (du,X)=0$. Thus, $Ker (d_*)=Im(m_*)$.

2. $Ker(i_*)=Im (d_*)$ because $i_*$ is the natural map of cocycles onto cohomology classes
in~$H^*_X$. For $u\in Ker(i_*)$, we have $u=dv, v\in\Lambda^k_X$. So we have $[u]=d_*[v]$
for the corresponding classes.

3. As $u\in Ker(j_*)$, it follows that $v=(u,X)=\nabla_X^k(w)=d(w,X)\pm (dw,X)$,
where $(w,X)=0$. So we have $v=(dw,X$. We see that $(u-dw,X)=0$. Hence, $u$ is
equivalent to an element from $ Im (i_*)$ in $H^{k+1}_X$.

4. Suppose that $v\in Ker(h_*)$ and $(u,X)=v, u\in\Lambda^{k+1}(M)$. This implies
that $du=dq, q\in\Lambda^{k+1}_X$. Replace $u$ by $u'=u-q$ with $du'=0$ .
Thus, $v=j_*(u')$, and hence $Ker(h_*)=Im(j_*)$.

5. Take $a\in Ker(g_*)$. This means that $da=0$ and $a=db$. We see
that $a=h_*(v)$ for $v=(b,X)$.

6. The last map $H^{k+2}_X\rightarrow H^{k+2}(M)$ is obviously an epimorphism.

\smallskip
\noindent
The theorem is proved.

\smallskip
\noindent
{\it Proof of Proposition\/~{\rm 1} above.} Our exact sequence for $k=n-1$
yields an isomorphism
$$\{0=H^n(M/X)\}\rightarrow \{H^n_X= Coker (\nabla^n_X)\}\rightarrow
C^{n-1}_ X\rightarrow \{H^{n+2}_X=0\}.$$
Taking into account the isomorphism $Coker(\nabla^{n}_X)= Coker(\nabla^0_X)$
for a measure-preserving flow with an invariant $C^{\infty}$ volume
form $\sigma$, we see that our assertion is proved.

\section{Riemann Surfaces and Hamiltonian Systems}

Consider now a compact nonsingular complex algebraic curve (Riemann surface) $M$
and any generic holomorphic one-form $\omega=\omega'+i\omega''$ on it.
The equation $dH=\omega'=0$  defines a Hamiltonian foliation.
Any smooth nonzero 2-form $\Omega$ defines a flow with vector field $X$ and
multivalued Hamiltonian $H$ such that $(\omega',X)=(dH,X)=0$.
Therefore $\nabla^{1,*}_X(\omega')=0$ and $\omega'\in \Lambda^1(M/X)$ where $d\omega'=0$.
It is easy to see that no other closed invariant one-form exists if the genus $g$ of $M$ is~$\geq 2$.
For $g=1$, there is exactly one more real closed invariant constant one-form $\omega''$ but it does
not belong to the space $\Lambda^1(M/X)$. In the generic case, the forms $\omega',\omega''$
have exactly $2g-2$ simple zeroes defining the saddle points of our system, and the system
is ergodic. Every function of the form $\nabla^0_X(f)=X^id_if$, $i=1,2,$ vanishes
at all saddle points as well as the space $Im(i_X)\subset Ker (i_X)=C^{\infty}(M)=\Lambda^0(M)$.
The space $C^0_X$ is infinite dimensional. We prove
\begin{cor}
There is the following exact sequence for a Hamiltonian system defined by
a generic Riemann surface $M$ and a holomorphic one-form $\omega$:
$$0\rightarrow \{R=H^1(M/X)\}\rightarrow
\{R^{2g}=H^1(M)\}\rightarrow C^0_X\rightarrow $$
$$\rightarrow \{H^2_X=Coker(\nabla^0_X)\rightarrow
\{R=H^2(M)\}\rightarrow 0.$$
The operator
$$h_*:C^0_X\rightarrow H^2_X(M)$$
is a Fredholm (``Noetherian'') operator with finite dimensional spaces
$Ker(h_*)$ and $Coker(h_*)$ whose dimensions are equal to $2g-1$
and $1$, respectively. The index of this operator is equal to
the Euler characteristic of the Riemann surface, that is, $-index(h_*)=2-2g$.
In particular, for $g\geq 2$ we have $\dim(Coker(\nabla^0_X))=\infty$ for this system
in the space of $C^{\infty}$ functions.
\end{cor}

\begin{rem}Even the case $g=1$ might present some difficulties.
For a generic slope $\alpha$, we need to invert the operator
$$A=\nabla^0_X=\partial_x +\alpha\partial_y:C^{\infty}(T^2)\rightarrow C^{\infty}(T^2),$$
where $x,y$ are defined modulo $1$ on the torus $T^2=M$. For generic non-Liouville numbers
$$|\alpha-m/n|>const_p(n)^{-p}$$ where $p>2$, we see that this operator
is invertible on the subspace consisting of $f\in V\subset C^{\infty}$
such that $\iint_{T^2}f(x,y)\,dx\,dy=0$. Indeed, the Fourier coefficients $a_{m,n}$
tend to zero with superpower speed as $m^2+n^2\rightarrow\{\infty\}$. For the $C^{\infty}$-function
$$h=Af=\sum_{)m,n)\neq 0} a_{m,n}\exp\{2\pi i (mx+ny)\},$$
we have small denominators for $f=A^{-1}(g)$, namely,
$$f=\sum 1/(2\pi i)a_{m,n} (mx+ny)^{-1}.$$
For a generic number $\alpha$, this series converges to a $C^{\infty}$-function $f$,
so we have  $Coker(\nabla^0_X)=\mathbb R$.
However, we may have $\dim(Coker(\nabla^0_X))=\infty$ for the exceptional numbers $\alpha$
with exponentially fast rational approximation.
\end{rem}

For generic flows $X$ of this kind, Forni~\cite{F} investigated the image of the operator $\nabla^0_X$
in the Sobolev spaces $H^{r-1}(M)$  for some natural number $r\geq 1$ (as he claims, one can take $r=5$).
For functions $f$ from that space (compactly supported outside of the singularities of $X$), he proved
the existence of a solution in the Sobolev space of distributions $H^{-r}_{loc}(M\setminus (X=0))$
with suitable estimates of the norm. For every number $s>2r-1$, there exist a finite number of $X$-invariant
distributions $D_q\in H^{-s}(M),q=1,...,n_s$ such that for every smooth function  $f\in H^s(M)$ compactly
supported in $M\setminus (X=0)$, satisfying the linear restrictions $D_q(f)=0$ for all $q=1,...,n_s$,
there exists a solution $u\in H^{s-2r+1}(M)$ with suitable estimates.
The number $n_s$ is not effective. Forni has also investigated the space of $X$-invariant
distributions in negative Sobolev spaces.

In the community of people working in dynamical systems,
the differential equation $\nabla^0_X(u)=f$ is called ``cohomological.''
Whatever it might mean, it has nothing to do with usual cohomology:
this terminology originated from the study of systems with discrete time
where the equation has a group-theoretic ``cohomological'' form $g(Tx)-g(x)=f$.

The work \cite{F} is based on some hard analysis.
Our arguments are very simple and general. They look as the simplest
natural analog of the de Rham cohomology theory for dynamical systems.
By definition, they work  for $C^{\infty}$-forms only.

To prove the corollary, let us denote the space $Coker(\nabla^0_X)$ by $W$ and rewrite
our exact sequence in the form
$$0\rightarrow \mathbb R\rightarrow \mathbb R^{2g}\rightarrow C^0_X\rightarrow W\rightarrow \mathbb R \rightarrow 0.$$
For $g=1$ and a non-Liouville slope $\alpha$, we have $\dim(W)=1$.
For $g\geq 2$, this exact sequence cannot be realized by finite-dimensional spaces because
the index of the operator $h_*:C^0_X\rightarrow W$ is strictly positive whereas $\dim (C^0_X)\leq \dim (W)$
because $C^0_X\subset W$. Our statement is proved.

\smallskip
\noindent
{\bf Explicit Construction:}

{\it Our exact sequence combined with the isomorphism $$v\rightarrow (v-\bar{v})\sigma$$
of the spaces $Coker(\nabla^0_X)/const$ and $Coker(\nabla^2_X)/\mathbb R$ plus the
knowledge of the subspace $Ker(g_*)\subset H^1_X(M)$ leads to an explicit construction
of the infinite-dimensional subspace $W'\subset W=Coker(\nabla^0_X)$.
First of all, we construct the subspace $\mathbb R^{2g-1}=V_0\subset W'$ as $V_0=j_*H^1(M)\subset W$. We have
$v_0=(u_0,X)\in V_0$ for the harmonic forms $u_0\in H^1(M)/(dH=\omega)'$
by the definition of the exact sequence. After that we take $v_0\sigma$ for
an invariant area form~$\sigma$. This gives a natural realization
of the isomorphism between $Coker(\nabla^0_X)$ and $H^2_X(M)$.
To get the projection on the subspace of exact $2$-forms,
we calculate the integral mean $\bar{v_0}=\iint_{M}v_0\sigma$ and subtract it
$$v_0\sigma\rightarrow (v_0-\bar{v_0})\sigma.$$
Next, solve the equation $$du_1=(v_0-\bar{v_0})\sigma.$$
By adding functions  $v_1=(u_1,X)$ to the subspace $V_0$, we obtain a larger subspace
$$ V_1= V_0+A(V_0),$$
where $A(v_0)=v_1$ was constructed above.
After that we apply to $v_1\in A(V_0)$ the same procedure
$$v_1\rightarrow v_2=(d^{-1}[(v_1-\bar{v_1}) \sigma)],X)$$
or $A(v_1)=(u_2,X)$ where $du_2=(v_1-\bar{v_1})\sigma$.
Iterating this construction, we obtain an infinite dimensional subspace
$$W'=V_0+A(V_0) +\ldots+A^m(V_0)+\ldots\subset W$$
because all subspaces $A^m(V_0)$ are linearly independent in $W=Coker(\nabla^0_X)$.
Each of them has dimension equal to the absolute value of the index of operator $g_*$,
which is $2g-2$.}

\section{Higher dimensional manifolds}

Consider now an orientable $n$-dimensional manifold $M$ with a
smooth flow (vector field) $X$. Using Theorem 1, we obtain the
following result:

\begin{cor}
For $k=1$, the exact sequence map $C^0_X\rightarrow H^2_X$ is a
Fredholm operator whose index is equal to $b_1(M)-|H^1(M/X)|-b_2(M)$.
The dimension of its kernel is equal to the difference $|H^1(M)|-|H^1(M/X)|$
where $|H^1(M)|=b_1(M)$, and the dimension of its cokernel is equal
to $|H^2(M)|=b_2(M)$.
\end{cor}

For $k\geq 1$, the exact sequence from Theorem 1 takes the  form
$$0\rightarrow H^{k+1}(M/X)\rightarrow H^{k+1}_X\rightarrow C^k_X\rightarrow H^{k+2}_X\rightarrow H^{k+2}(M)\rightarrow 0.$$
For measure-preserving flows, we also have $H^n_X=Coker(\nabla^0_X)=Coker(\nabla^n_X)$,
as it was pointed out before. This gives us the {\bf cokernel complex} $(C,d_C)$:
$$0\rightarrow C^0_X\rightarrow C_X^1\rightarrow \ldots
C^{n-1}_X \rightarrow 0$$ or $$0\rightarrow C^0\rightarrow \ldots
C^{n-1}\rightarrow 0.$$
The coboundary operator $d_C$ extracted from the exact sequence of Theorem~1 above,
can be expressed by the formula
$$d_C(v)=i_Xdi_X^{-1},\quad v\in C^j.$$

\begin{thm}
The cohomology groups $H^i_C$  of the complex $C$ can be described by the following
long exact sequence\/{\rm :}
$$\ldots\rightarrow H^{k+1}(M/X)\rightarrow H^{k+1}(M)\rightarrow H^k_C
\rightarrow H^{k+2}(M/X)\rightarrow H^{k+2}(M)\rightarrow\ldots$$
In particular, the maps $h_k:H^k_C\rightarrow H^{k+2}(M/X)$ are
Fredholm operators,
and their indices satisfy
$$\sum_{k\geq 0}(-1)^k {\mathrm{Index}}\, (h_k)=1-b_1(M/X)-\sum_{k\geq
0}(-1)^kb_k(M).$$
\end{thm}

In order to prove this theorem, let us define all the homomorphisms involved.
The map $H^j(M/X)\rightarrow H^j(M)$ is obvious because our complex $\Lambda^*(M/X)$
is a subcomplex in $\Lambda^*(M)$ with the same differential. The map
$H^{j+1}(M)\rightarrow H^j_C$ is defined by setting $u\rightarrow i_X(u)=v$
for any representative of a cohomology class. The map $H^j_C\rightarrow H^{j+2}(M/X)$
is defined by setting $v\rightarrow d i_X^{-1}(v)$. This is well-defined because
$C=Im(i_X)/Im(\nabla^*_X)$. The theorem follows now from the same elementary arguments
as the standard exact homological sequences in basic algebraic topology.

\begin{ex}
Consider a level set $M_c:|p|=c>0$ in a cotangent bundle $T_*(N)$.
The geodesic flows on the manifold $T_*(N)$ are described by
the Hamiltonians $H=1/2|p|^2$ corresponding to Riemannian metrics.
They are homogeneous functions of order two in the momentum variables.
So in addition to the standard 2-form $\sum dp_i\wedge dx^i|_{H=c}=\Omega_c\in
\Lambda^2(M_c/X)$, we  also have the corresponding standard 1-form
$\sum p_idx^i|_{H=c}=\omega\in \Lambda^1_{inv}$ where $d\omega=\Omega_c$.
However, $i_X(\omega)\neq 0$. Assuming that generically there are no other
smooth forms in the complex $\Lambda^*(M_c/X)$, we have $H^{2j}(M_c/X)=\mathbb R,j=1,...,n-1,$
but $H^{2j}_{inv}(M_c)=0$. However, nobody has proved this statement yet.
Performing a small potential perturbation $H'=1/2|p|^2+U(x)$, we remove
the forms $\omega\wedge\Omega^j, j<n-1,$ out of the space $\Lambda^1_{inv}(M_c)$.
So the nontrivial invariant cohomology appears in the same even dimensions for
the perturbed field $X'$, i.e., $H^{2j}_{inv}(M_c)=\mathbb R$. For both fields
$X,X'$ we have $H^{2n-1}_{inv}(M_c)=\mathbb R$, where $ n=\dim N$.
At the moment, this ``generic'' statement is non-rigorous.
\end{ex}

\begin{ex}
Let $M_c\subset M$ be a compact energy level $H=c>0$ in a symplectic
$2n$-manifold. For generic Hamiltonians, the complex $\Lambda^*(M_c/X)$ consists
of the powers $\Omega_c^j,j=0,1,\ldots, n-1$, so we have $H^{2j}(M_c/X)=\mathbb R$.
It is worth mentioning that the case of geodesics is non-generic, and an
additional invariant one-form appear in the complex $\Lambda^*_{inv}$.
\end{ex}

\begin{ex}
Consider a compact symplectic $2n$-manifold $(M,\Omega)$ such that $H^1(M)\neq 0$.
Every closed 1-form $\omega=dH$ defines a Hamiltonian flow with energy-foliation $\omega=0$
or $H=const$. Here $H$ is a multivalued function on $M$ which becomes single-valued on some
abelian covering $p:\hat{M}\rightarrow M$, $p^*(\omega)=dH$. Generically, we have $H^1(M/X)=\mathbb R$,
generated by the class of~$\omega$. For $2$-manifolds, this example was already discussed
in Section~2 above. All forms $\Omega^j$ and $\omega\wedge\Omega^j$ lie in the invariant
subcomplex $\Lambda^*_{inv}$. Let us prove the following result.

\begin{lem}
The forms $\mu_j=\omega\wedge\Omega^j$ lie in the subcomplex
$\Lambda^*(M/X)$, i.e., they are invariant and such that
$(\mu_j,X)=i_X(\mu_j)=0$.
\end{lem}

\noindent
{\it Proof.} In appropriate local coordinates, we have
$\Omega=\sum_{i=1}^{n-1} dx^i\wedge dp^i+dc\wedge dt$ and
$H=c,\omega=dc$. So we conclude that $ dc\wedge \Omega^j$
does not contain differential ~$dt$. The lemma is proved.

\smallskip
In the generic case, there are no additional forms in $\Lambda^*(M/X)$,
so we conclude that generically the odd-dimensional homology groups $H^{odd}(M/X)$ are non-trivial.
\end{ex}

\begin{ex}
Let $M=T_1(M^2_g)$ be a space of unit linear elements on a closed
surface of genus $g$ with a metric of negative curvature. Its geodesic
flow $X$ is a strongly ergodic system (even a K-system).
It was proved many years ago by Livshitz (see \cite{L}) that the equation
$$\nabla^0_X(f)=g\in C^{\infty}(M)$$
is solvable for functions $g$ whose integrals along closed geodesics are zero.
The solution $f$ belongs to the space $C^1(M)$ as it was proved in \cite{L}.
This result is true for all Anosov systems. For constant negative curvature
surfaces $N$ and  $M=T^1(N)$, Kazhdan and Guillemin proved in \cite{KG}
that the solution is in~$C^{\infty}$. The variable curvature case was settled
in~\cite{CEG} and~\cite{LMM}. So the space $C_X^0=C^{\infty}(M)/Im(\nabla^0_X)$
is infinite-dimensional. Every closed geodesic gives a nontrivial linear functional
$$I_{\gamma}:f\rightarrow \oint_{\gamma}fdl, f\in C^{\infty}(M)/Im(\nabla^0_X)$$
on the space $C^0_X$, and so we have infinite number of linearly independent functionals.
How can one describe the topology of this factor-space? How does it depend on the discrete group?
The answers to these questions are unknown.
\end{ex}

\noindent
{\bf Problem:} Prove rigorously that generically (i.e., after a generic perturbation)
there remains no invariant differential form except those described in the above examples.

\section{Exotic Poisson Manifolds. Dynamical Systems on $PSL_2(R)/\Gamma$.}

The most general class of $C^{\infty}$ Hamiltonian systems is described by the
following definition. Let $M$ be a $C^{\infty}$-manifold equipped with a
Poisson tensor field. We call it a Poisson manifold. Locally the Poisson
tensor can be written in the form $h^{ij}(x)$. The corresponding ``Poisson bracket''
of smooth function is defined locally by $\{f,g\}=h^{ij}(x)f_{x_i}g_{x_j}$
and satisfies the following  relations:
$$\{f,g\}=-\{g,f\},\;\{f,gh\}=h\{f,g\}+g\{f,h\},\;\{\{f,g\},h\}+(cyclic)=0$$
In particular, we have $h^{ij}=\{x^i,x^j\}$ for the local coordinates $x^i$.
In the non-degenerate (``symplectic'') case, these conditions imply that
the 2-form $\Omega=\sum h_{ij}(x)dx^i\wedge dx^j$, where $h_{ij}h^{jk}=\delta^k_i$,
is closed.

The simplest natural classes of Hamiltonian systems correspond to single-valued
Hamiltonian functions $H$ by the formula $$df/dt=\{f,H\}$$ for a function $f\in C^{\infty}(M)$.
More general Hamiltonian systems are given by closed 1-forms $dH$ by the same local
formula $$df/dt=\{f,H\},$$ where $H$ is locally well-defined modulo an additive constant
(which does not affect the equation). No other Hamiltonian system exists on a symplectic
manifold $M$ with non-degenerate Poisson tensor.

However, this is no longer true for more general Poisson manifolds with degenerate Poisson tensor.
The Poisson analog of the Darboux theorem claims (see \cite{NT}) that for every point
$x\in M$ there exists an open set $U\ni x$ and a local coordinate system
$(x^1,\ldots,x^{2n},y^1,\ldots,y^k)$  in~$U$ such that
$$\{x^i,x^{n+i}\}=1, \{y^p,y^q\}=h^{pq}(y),$$
where $h^{pq}(0)=0$, and the other Poisson brackets of coordinates are zero.
Here the number $2n$ is equal to the rank of the Poisson matrix at the origin~$x^i=0,y^p=0$.
Let $$h^{pq}(y)=C^{pq}_ry^r+O(|y|^2.$$ It is easy to prove the following

\begin{lem}
The finite-dimensional Lie algebra $L$ defined by the structural constants $C^{pq}_r$ with respect
to a basis $e^1,\ldots,e^k$ so that
$$[e^p,e^q]=C^{pq}_re^r$$
is an invariant of the Poisson tensor at the point $x\in M$.
\end{lem}

This Lie algebra invariant sometimes determines the Poisson structure completely,
and the higher terms are inessential. This problem was studied in~\cite{DZ}.
These results are especially effective for low dimensional Lie algebras.

A lot of people have studied  special ``Lie-Poisson'' structures for which
the Poisson tensor linearly depends on the coordinates. In this case $M=L^*=\mathbb R^k$,
and the basis elements $e^i$ can be viewed as linear functions on the manifold $M$.
In particular, there exists a complete set of ``Annihilators'' or ``Casimirs''
which are single-valued functions $f_1,\ldots,f_t$ (homogeneous polynomials
in the variables $e^1=y^1,\ldots,e^k=y^k$) such that $\{f_m,y^p\}=0$
for all $p=1,2,\ldots,k$. Any Casimir (even locally defined in some open set in $M$)
is functionally dependent on this family $(f_m)$. Casimirs define trivial Hamiltonian
systems. They are  integrals of motion  for all hamiltonian systems in $M$.
Their common levels $f_m=c_m$  define the so-called ``Symplectic Foliation''
of the Poisson manifold $M$.

What  do we have for more complicated Poisson manifolds? How do Casimirs look like?
In the simplest cases they are either one-valued functions (as in the Lie-Poisson manifolds $M=L^*$ above)
or closed one-forms. The last case was realized in  Quantum Solid State Physics. ``Semiclassical''
electrons move under the influence of a magnetic field in the space of quasimomenta~$p\in M$.
The latter is the Euclidean space $\mathbb R^3$ factorized by the reciprocal lattice, i.e., $M=T^3$.
Its points numerate the quantum states of an electron when a branch of the energy dispersion relation
is fixed and the magnetic field is absent.
Switching on a constant magnetic field, we get a  Poisson tensor $B=(B^{pq})=const$.
It has a multivalued Casimir (i.e., a constant one-form $*B$ on $T^3$ dual to the  magnetic field $B$).
As it is written in the  physics literature, electrons move in the plane orthogonal to the magnetic
field regarded as a vector~$*B$. The Hamiltonian $\epsilon(p):T^3\rightarrow R$ is a single-valued
Morse function in this case (it is the energy dispersion relation extracted from quantum theory).
So the trajectories are the sections of the ``Fermi-surface'' $\{\epsilon(p)=\epsilon_F\}\subset T^3$
by the planes orthogonal to the magnetic field.

As A.Ya.Maltsev has pointed out to me, Casimirs  might be more complicated than closed one-forms
in general. Let us now consider the special class of Poisson tensors of constant rank~$2n$.

\begin{lem}
Any Poisson tensor defines a completely integrable orientable foliation on the manifold $M$,
which we call the ``Casimir''  or ``annihilator'' foliation, and a ``symplectic'' 2-form
$\Omega\in \Lambda^2(M)$ such that the restriction of $\Omega$ on any leaf is non-degenerate and closed.
This foliation is defined by the ideal $I\subset\Lambda^*(M)$ such that $dI\subset I$.
It is generated by the collection of one-forms $\omega_j\in I$ equal to zero on every leaf $\omega_j=0$.
We have $d\Omega\in I$ and $d\omega_j\in I$, i.e., $dI\subset I$ . Vice versa, every completely integrable
foliation equipped with a 2-form $\Omega$ which is closed and non-degenerate on every leaf, defines a
Poisson tensor. Two different 2-forms $\Omega'$ and $\Omega$ define the same Poisson tensor
if $\Omega'-\Omega\in I$.
\end{lem}

\begin{cor}Every orientable foliation with two-dimensional leaves defines at least one Poisson tensor.
In this case, any two Poisson tensors with the same Casimir foliation are proportional to each other
so that $h^{ij}_1(x)=h^{ij}_2(x)f(x)$ with $f(x)\neq 0$. Therefore, the corresponding Hamiltonian systems
differ only by a time change and have the same trajectories.
\end{cor}

To prove this statement, we use a Riemannian metric. It defines an area element along the leaves
which can be regarded as a 2-form $\Omega\in \Lambda^2(M)$ because everything is orientable.
The restriction of this  2-form to the leaves is unique up to multiplication $\Omega\rightarrow g(x)\Omega$,
where $g>0$.

Let us discuss several important examples.

1. In a series of papers started in 1982 we investigated the
following class of dynamical systems (see the survey in \cite{ND}):

A constant skew-symmetric $n\times n$-matrix $B$ of rank~$2$ defines
a Poisson structure on the $n$-torus $T^n=R^n/\Gamma$, where
$\Gamma$ is an Euclidean  lattice of rank~$n$. The kernel of the
Poisson matrix consists of $n-2$ linear functions $l_1,...,l_{n-2}$.
A Hamiltonian $H:T^n\rightarrow\mathbb R$ defines a system whose
trajectories are intersections of the level sets $H=c$ with the
$2$-plane $l_1-a_1,...,l_{n-2}=a_{n-2}$. One might say that every
trajectory is a level set of the quasiperiodic function $H$ on the
plane. This quasiperiodic function has $n$ frequencies
(quasiperiods). {\bf What can one say about the topology of the
levels of quasiperiodic functions on the plane?} This problem
appears in several different areas.

For $n=3$ this problem appeared from the solid state physics
describing the motion of electrons along the Fermi surface in the
single crystal normal metal under the influence of a strong magnetic
field~$B$. It should be added ``semiclassically'' to this purely
quantum picture. For low temperature this approximation is valid for
magnetic fields in the interval $1 te<|B|<100te$. There exists
complicated examples of dynamical systems on the Fermi surface in
this problem but they are nongeneric: generically every open
trajectory lives in some piece of genus one. Some remarkable
observable integer-valued topological characteristics of electrical
conductivity in  magnetic field were found here for the generic case
(see \cite{ND}; A.\,Zorich, I.\,Dynnikov, S.\,Tsarev, A.\,Maltsev
and R.\,Deleo have worked in this area jointly with the present
author since 1980s). Similar features were recently found also for
$n=4$ in~\cite{ND}. Further applications of this problem to the
physics of quantum surfaces were found by A.\,Maltsev for all $n$
(see the references in~\cite{ND}).

\smallskip
2. Let $M=T_1(M^2_g)$ be the space of unit linear elements on a
compact surface of genus g with a Riemannian metric of constant
negative curvature. We have $M=PSL_2(R)/\Gamma$, where $\Gamma$
is a cocompact discrete group acting from the right. The space
of right-invariant differential forms is generated by the $1$-forms
$\omega_0,\omega_{\pm}$ such that
$$d\omega_0=\omega_+\wedge\omega_-,d\omega_{\pm}=\pm\omega_0\wedge\omega_{\pm}.$$
This forms are dual to the right-invariant vector fields
$e_0.e_{\pm}$ such that
$$[e_+,e_-]=e_0, [e_0,e_{\pm}]=\pm e_{\pm}$$
The {\bf important Poisson tensors} $h_{\pm}$ are given by the formulas
$$h_+=e_0\wedge e_-,\;h_-=e_0\wedge e_+,$$
and their symplectic foliations by $\omega_+=0$ and $\omega_-=0$, respectively.

The most famous dynamical systems here are the following:

\smallskip
a. The geodesic flow. Its trajectories are given by the intersection of two completely
integrable ``Anosov'' foliations $$\omega_+=0,\;\omega_-=0.$$

b. The horocycle-$(+)$ flow whose trajectories are given by the
equations
$$\omega_+=0,\;\omega_0=0.$$

c. The horocycle-$(-)$  flow whose trajectories  are given by the equations
$$\omega_-=0,\;\omega_0=0$$

\smallskip
Let us point out that the form $\omega_0$ is closed on the leaves
$\omega_{\pm}=0$. Besides, the restriction of the 1-form
$\omega_0$ to the geodesics contained in these leaves is
proportional to the element of length with a non-zero constant
factor. Therefore, we have proved the following result.

\begin{pr}
The ($\pm$)-horocycle flows are Hamiltonian systems with the same
Hamiltonian 1-form $\omega_0$ but different Poisson tensors with the
the (Anosov) Casimir foliations $\omega_{\pm}=0$, respectively. For
the topologically non-trivial (cylindrical) leaves containing the
closed geodesics $\gamma$, the period of  hamiltonian 1-form is
always nonzero and equal (after the proper choice of the dimensional
constant)  to the length of the geodesic
$$\oint_{\gamma}\omega_0=l(\gamma).$$
\end{pr}

The whole space of differential forms on $M$ can be described in
this basis. Let us denote the ring of $C^{\infty}$-functions
$\Lambda^0(M)$ by $S$. The space of differential forms is then
represented as an external power of the 3-dimensional free
$S$-module
$$\Lambda^*(M)=\Lambda^*[S_3],$$
where $S_3$ is linearly generated over the ring $S$ by the forms
$\omega_0,\omega_{\pm}$.

a. For the geodesic flow $X$ we know that the image of the operator $\nabla^0_X$
consists of the functions $g$ such that $$\oint_{\gamma}g\,dl=0$$ for all closed geodesics~$\gamma$.
For the right-invariant forms, we have
$$\nabla^1_X(\omega_{\pm})=\pm\omega_{\pm},\nabla^1_X(\omega_0)=0,\;
\nabla^3_X(\omega_0\wedge\omega_+\wedge\omega_-)=0,$$
$$i_X(\omega_0)=1,\; i_X(\omega_{\pm})=0,\; i_X(\omega_0\wedge
\omega_+\wedge\omega_-)=\omega_+\wedge\omega_-,$$
$$i_X(\omega_0\wedge\omega_{\pm})=\pm\omega_{\pm},\; i_X(\omega_+\wedge\omega_-)=0.$$
Taking into account the identity (no signs)
$$\nabla_X(a\wedge b)=\nabla_X(a)\wedge b+a\wedge \nabla_X(b),$$
we can calculate this operator for all forms. In particular, we get
$$\nabla^2_X(\omega_+\wedge\omega_-)=0.$$
This equality reflects the fact that the geodesic flow is a
restriction of a 4D Hamiltonian system to the energy level. This
system is non-Hamiltonian for the 3D Poisson structure  described
above. Even the trajectory foliation $\omega_+=0,\omega_-=0$ does
not correspond to any Hamiltonian system. The reason for this is the
following: geodesics are located on the leaves $\omega_+=0$ but the
form $\omega_-$ restricted to these leaves is non-closed because
$d\omega_-=\omega_0\wedge \omega_-$. In order to construct
hamiltonian,  one needs to find a global 1-form $f\omega_-$ such
that $f\neq 0$ and $d(f\omega_-)=0$ on every leaf. However, this is
impossible because there are limit cycles on the cylindrical leaves
corresponding to the closed geodesics. The set of all closed
geodesics is everywhere dense.

Consider the operator $A=\nabla_X^0$ acting on $S=C^{\infty}(M)$.
The equation
$$\nabla_X^1(f\omega_0+g\omega_++h\omega_-)=A(f)\omega_0+(g+A(g)) \omega_++(h-A(h)\omega_-=0$$
is solvable only for $f=const, g=h=0$ because $\pm 1$ cannot be realized as eigenvalues
of the skew-symmetric operator $A$. So we have $$\Lambda^1_{inv}=\mathbb R,\; \Lambda^1(M/X)=0,$$
and $\Lambda^0_{inv}=\Lambda^0(M/X)=\mathbb R$ because the flow is ergodic. Our arguments imply
that
$$\Lambda^2_{inv}=\Lambda^2(M/X)=\mathbb R\;\text{ and }\;\Lambda^3_{inv}=\mathbb R.$$
Thus, the exotic homology
$$H^j_{inv}=\mathbb R, j=0,3,\; H^j_{inv}=0,j=1,2,$$
$$H^j(M/X)=\mathbb R, j=0,2,\; H^1(M/X)=0.$$
Using the exact sequence from Theorem~2, we can compute the homology of
the cokernel complex $C$:
$$H^2_C=H^3(M)=\mathbb R$$
$$0\rightarrow\mathbb R^{2g}\rightarrow H^0_C\rightarrow H^2(M/X)\rightarrow
\mathbb R^{2g}\rightarrow H^1_C\rightarrow 0,$$
where
$$H^2(M/X)=\mathbb R, H^1(M)=H^2(M)=\mathbb R^{2g}.$$

b. Similar calculations can be performed for the $\pm$-horocycle flows
with vector fields $X_{\pm}$. Let us use the notation $i_{\pm}, \nabla_{\pm}$
for the corresponding operators $i_{X_{\pm}},\nabla_{X_{\pm}}$.
Then by the definition of horocycles, we have
$$i_+(\omega_0)=i_+(\omega_+)=0,i_+(\omega_-)=1,$$
$$i_+(\omega_0\wedge\omega_+)=0,i_+(\omega_0\wedge\omega_-)=\pm\omega_0,
 i_+(\omega_+\wedge\omega_-)=\pm\omega_+.$$
This implies that
$$\nabla^1_+(\omega_0)=\omega_+,\nabla^1_+(\omega_)=\pm\omega_0,\nabla^1_+(\omega_+)=0.$$
In particular, $$\nabla^2_+(\omega_0\wedge\omega_+)=0$$
because this system is obtained by the restriction of the Hamiltonian system on the $4D$
manifold $T^*(M_2)$ to the energy level $T_1$. This Hamiltonian system has the same
Hamiltonian function $H=1/2|p|^2$ as the geodesic flow, but the symplectic structure
is changed, as should always be done after switching on a magnetic field. (In order
to get the correct symplectic form, the magnetic field should be added to the original
symplectic form.) In our situation, the magnetic field is equal to the Gaussian curvature
form on the surface lifted to the phase space $T^*(M_2)$.

For the ``symplectic'' differential 2-form $\Omega'=\omega_0\wedge\omega_-$ associated
with the Poisson tensor on the 3-manifold $M$, we have $\nabla^2_+(\Omega')=\omega_0\wedge\omega_+$,
which is $0$ on the leaves $\omega_+=0$.

We conclude also that
$$\omega_+\in \Lambda^1(M/X_+)\subset   \Lambda^1_{inv}$$
and
$$\omega_0\wedge\omega_+\in \Lambda^2(M/X_+)\subset
 \Lambda^2_{inv}, \omega_0\wedge\omega_+\wedge\omega_-\in \Lambda^3_{inv}.$$
Taking into account these formulas combined with the facts that the
flow is ergodic  and the operator $\nabla_+$ is skew-symmetric, we
can deduce that all $X_+$-invariant forms are right-invariant. So we
have
$$H^0(M/X_+)=R,H^j(M/X_+)=0,j=1,2,$$
$$H^j_{inv}=R,j=0,3,H^1_{inv}=H^2_{inv}=0.$$
Using  the Theorem 2, we obtain the following result for the ``cokernel cohomology'':
$$H^j_C=H^{j+1}(M),j=0,1,2.$$
For all $k>0$, the operators $\nabla^k_+$ on the subspaces of right-invariant forms
are nilpotent
with a single Jordan cell only.
For $k=0$, the flow has no periodic orbits and is ergodic.
The image $Im(\nabla^0_+)\subset C^{\infty}(M)$ is certainly
contained in the standard subspace of functions with zero integral
$\int\!\cdots\!\int_Mfd\sigma=0$, where the volume element
$d\sigma=\omega_0\wedge\omega_+\wedge\omega_-$. Flaminio and Forni
proved in~\cite{FF} that the image $Im(\nabla^0_+)$ has infinite
codimension and is closed in~$C^{\infty}(M)$. The same result is
valid  for the second horocycle flow~$X_-$.

For any negatively curved compact  2-manifold $M^2_g$ of
non-constant curvature, we also have a pair of foliations
$\omega_{\pm}$ on $M=T_1(M^2_g)$ and a natural definition of the
horocycle flows. These foliations are generically such that every
leaf is smooth but the whole family is only  $C^1$. So our analogs
of the horocycle flows and Poisson Tensors are $C^1$ only. The form
$\omega_0$ needed here is naturally defined by the geodesic flow as
the restriction of the canonical one-form $\sum p_jdx^j$ from
$T_*(M^2_g)$ to the level $H=1$. All previous results are valid
here.

Interesting results can be obtained for Hamiltonian systems
with Anosov-type Casimir foliations and single-valued Hamiltonians
$H:M\rightarrow\mathbb R$. They lead to ``locally Hamiltonian'' dynamical systems
on the surfaces $H=c$ (i.e., such that their singular points are like in Hamiltonian
systems). These systems may have only very special limit cycles that are non-homotopic to zero
in $M$ and have lengths bounded from below. There are very interesting problems here,
which we postpone to future publications.

\end{document}